\documentclass{article}
\usepackage{graphicx}
\usepackage{amsmath}
\usepackage{amsfonts}
\usepackage{amssymb}
\newtheorem{theorem}{Theorem}

\newtheorem{lemma}[theorem]{Lemma}

\begin{document}

\title{On an alternate proof of Hamilton's matrix Harnack inequality for the Ricci flow}
\author{Bennett\ Chow\thanks{Research partially supported by NSF grant DMS-9971891.}\\Department of Mathematics\\University of California at San Diego}
\date{August 15, 2001}
\maketitle

In \cite{LY} a differential Harnack inequality was proved for solutions to the
heat equation on a Riemannian manifold. Inspired by this result, Hamilton
first proved trace and matrix Harnack inequalities for the Ricci flow on
compact surfaces \cite{H0} and then vastly generalized his own result to all
higher dimensions for complete solutions of the Ricci flow with nonnegative
curvature operator \cite{H2}. Soon afterwards, a matrix Harnack inequality for
the K\"{a}hler-Ricci flow under the assumption of nonnegative bisectional
curvature was proved by Huai-Dong Cao \cite{C}. In this paper, following a
suggestion of Richard Hamilton, we give an alternate proof of the matrix
Harnack inequality for the Ricci flow originally proved by him.\footnote{Yet
another proof is given in \cite{CC}. Hamilton's motivation for considering the
alternate proof given in the present paper is that it may possibly falicitate
extending his Harnack inequality to the variable signed curvature case,
especially in dimension three. See the end of \cite{H3} for some reasons for
studying the problem of extending the Harnack inequality.} In \cite{H2},
Hamilton proved that if $\left(  M^{n},g\left(  t\right)  \right)  $ is a
complete solution to the Ricci flow with nonnegative curvature operator, then%
\begin{equation}
\frac{\partial R}{\partial t}+\frac{R}{t}+2\nabla_{i}RV^{i}+2R_{ij}V^{i}%
V^{j}\geq0\label{trace-inequality1}%
\end{equation}
for any vector field $V.$ In particular, if in addition the Ricci curvature is
positive, then%
\[
\frac{\partial R}{\partial t}+\frac{R}{t}-\frac{1}{2}\left(  Rc^{-1}\right)
^{ij}\nabla_{i}R\nabla_{j}R\geq0.
\]
This trace Harnack inequality is a corollary of Hamilton's matrix Harnack
inequality, which we now recall (we follow Hamilton's sign convention for the
curvature tensor so that $R_{ijkl}=g_{kh}R_{ijl}^{h}$). Define tensors $P$ and
$M$ by%
\begin{align*}
P_{abc} &  \doteqdot D_{a}R_{bc}-D_{b}R_{ac}\\
M_{ab} &  \doteqdot\Delta R_{ab}-\frac{1}{2}D_{a}D_{b}R+2R_{acbd}R_{cd}%
-R_{ac}R_{bc}+\frac{1}{2t}R_{ab}%
\end{align*}
and consider the Harnack quadratic%
\[
Z\left(  U,W\right)  \doteqdot M_{ab}W_{a}W_{b}+2P_{abc}U_{ab}W_{c}%
+R_{abcd}U_{ab}U_{cd},
\]
where $U$ is any $2$-form and $W$ is any $1$-form. Hamilton's matrix Harnack
estimate says that if $\left(  M^{n},g\left(  t\right)  \right)  $ is a
complete solution to the Ricci flow with nonnegative curvature operator, then%
\[
Z\left(  U,W\right)  \geq0
\]
for all $U$ and $W.$ Taking $U=V\wedge W$ and summing over $W$ in an
orthonormal frame, one obtains (\ref{trace-inequality1}).

The idea of giving an alternate proof of it is as follows. Given any vector
$W,$ consider the $2$-form $U$ minimizing the quadratic form. This leads to
the consideration of a symmetric $2$-tensor $Z_{ab}$ which has a nice
evolution equation. In particular, $Z_{ab}$ is essentially a supersolution to
the heat equation so that the maximum principle implies that it is positive
definite, which is \emph{equivalent} to the matrix Harnack inequality. Tracing
this $2$-tensor yields a trace Harnack inequality: $Z\doteqdot Z_{aa}>0.$

Now we proceed with the details of the proof. Assume that the curvature
operator is positive and let $S_{abcd}$ denote the inverse of $R_{abcd}$ so
that%
\[
\sum_{e,f}S_{abef}R_{efcd}=I_{abcd},
\]
where $I:\wedge^{2}M\rightarrow\wedge^{2}M$ is the identity. First we rewrite
the matrix Harnack quadratic by completing the square:
\begin{align*}
Z\left(  U,W\right)   &  =M_{ab}W_{a}W_{b}-S_{ijkl}P_{ijp}W_{p}P_{klq}W_{q}\\
&  +R_{abcd}\left(  U_{ab}+S_{abij}P_{ijp}W_{p}\right)  \left(  U_{cd}%
+S_{cdkl}P_{klq}W_{q}\right)  .
\end{align*}
Since $Rm>0,$ the above expression is minimized when%
\[
U_{ab}+S_{abij}P_{ijp}W_{p}=0,
\]
in which case%
\[
Z\left(  U,W\right)  =\left(  M_{ab}-S_{ijkl}P_{ija}P_{klb}\right)  W_{a}%
W_{b}.
\]
Let
\[
Z_{ab}\doteqdot M_{ab}-S_{ijkl}P_{ija}P_{klb},
\]
which is a symmetric $2$-tensor.

The main computation is that in an orthonormal moving frame, the evolution of
$Z_{ab}$ is given as follows. Following the idea in \cite{H2}, as long as it
is nonnegative-definite, we may write the Harnack quadratic as a sum of
squares:%
\begin{align*}
R_{abcd} &  =\sum_{N}Y_{ab}^{N}Y_{cd}^{N}\\
P_{abc} &  =\sum_{N}Y_{ab}^{N}X_{c}^{N}\\
M_{ab} &  =\sum_{N}X_{a}^{N}X_{b}^{N}.
\end{align*}
Here, given a point $x\in M,$ $\left\{  Y^{N}\oplus X^{N}\right\}  _{N=1}^{m}$
is some orthogonal set of vectors in $\wedge^{2}M_{x}\oplus$ $\wedge^{1}M_{x}$
and $m=\dim\wedge^{2}M_{x}+\dim$ $\wedge^{1}M_{x}=\frac{n\left(  n+1\right)
}{2}.\medskip$

\noindent\textbf{Main Theorem. }\emph{If }$Rm>0,$\emph{ then the Harnack
quantity }$Z_{ab}$\emph{ satisfies the evolution equation:}%
\begin{equation}
\left(  \frac{\partial}{\partial t}-\Delta\right)  Z_{ab}=2S_{mntu}%
K_{avtu}K_{bvmn}+L_{a}^{NM}L_{b}^{NM}-\frac{2}{t}Z_{ab},\label{MT}%
\end{equation}
\emph{where}%
\begin{align*}
K_{avtu} &  \doteqdot S_{ijrs}P_{ija}\nabla_{v}R_{rstu}-\nabla_{v}%
P_{tua}+R_{tuaw}R_{vw}+\frac{1}{2t}R_{tuav}\\
L_{a}^{NM} &  \doteqdot2Y_{id}^{N}Y_{jd}^{M}S_{ijkl}P_{kla}+Y_{ah}^{N}%
X_{h}^{M}\newline -Y_{ah}^{M}X_{h}^{N}.
\end{align*}
\smallskip

Note that applying the maximum principle for symmetric 2-tensors in section 7
of \cite{H1} yields $Z_{ab}>0,$ which is equivalent to Hamilton's matrix
Harnack inequality. Here we used the fact that if $Rm>0,$ then as long as
$Z_{ab}$ is nonnegative definite, $L_{a}^{NM}$ is well-defined.

The bulk of the rest of this paper is devoted to proving the theorem. First
recall Hamilton's computations for the evolution equations of the Harnack
quadratic \cite{H2}.

\begin{lemma}
In an evolving local orthonormal frame field, we have

\begin{enumerate}
\item
\begin{equation}
(D_{t}-\triangle)R_{abcd}=2(B_{abcd}-B_{abdc}+B_{acbd}-B_{adbc})\,,
\label{Rabcd-evolution}%
\end{equation}
where $B_{abcd}=R_{aebf}R_{cedf}\,.$

\item
\begin{equation}
(D_{t}-\triangle)P_{abc}=-2R_{de}D_{d}R_{abce}+2R_{adbe}P_{dec}+2R_{adce}%
P_{dbe}\newline +2R_{bdce}P_{ade}\,. \label{Pabc-evolution}%
\end{equation}

\item
\begin{align}
(D_{t}-\triangle)M_{ab}  &  =\newline 2R_{cd}[D_{c}P_{dab}+D_{c}%
P_{dba}]\newline +2R_{acbd}M_{cd}\label{Mab-evolution}\\
&  +2P_{acd}P_{bcd}-4P_{acd}P_{bdc}\newline +2R_{cd}R_{ce}R_{adbe}-{\frac
{1}{2t^{2}}}\,R_{ab}\,.\nonumber
\end{align}
\end{enumerate}
\end{lemma}

We use this lemma in the computations below. We start by computing:%
\begin{align*}
\left(  \frac{\partial}{\partial t}-\Delta\right)  Z_{ab}  &  =\left(
\frac{\partial}{\partial t}-\Delta\right)  M_{ab}-\left(  \frac{\partial
}{\partial t}-\Delta\right)  S_{ijkl}P_{ija}P_{klb}\\
&  -S_{ijkl}\left(  \frac{\partial}{\partial t}-\Delta\right)  P_{ija}%
P_{klb}-S_{ijkl}P_{ija}\left(  \frac{\partial}{\partial t}-\Delta\right)
P_{klb}\\
&  +2S_{ijkl}\nabla_{p}P_{ija}\nabla_{p}P_{klb}+2\nabla_{p}S_{ijkl}\nabla
_{p}P_{ija}P_{klb}+2\nabla_{p}S_{ijkl}P_{ija}\nabla_{p}P_{klb}.
\end{align*}
Thus, besides the evolution equations for $M_{ab}$ and $P_{abc},$ we need the
evolution equation for $S_{ijkl},$ which is given by:
\begin{align*}
\left(  \frac{\partial}{\partial t}-\Delta\right)  S_{ijkl}  &  =-S_{ijmn}%
S_{klpq}\left(  \frac{\partial}{\partial t}-\Delta\right)  R_{mnpq}\\
&  -\left(  S_{ijrs}S_{mntu}S_{klpq}+S_{ijmn}S_{klrs}S_{pqtu}\right)
\nabla_{v}R_{rstu}\nabla_{v}R_{mnpq}\\
&  =-2S_{ijmn}S_{klpq}(B_{mnpq}-B_{mnqp}+B_{mpnq}-B_{mqnp})\\
&  -\left(  S_{ijrs}S_{mntu}S_{klpq}+S_{ijmn}S_{klrs}S_{pqtu}\right)
\nabla_{v}R_{rstu}\nabla_{v}R_{mnpq}.
\end{align*}
Substituting in those evolution equations, we obtain:%
\begin{align*}
\left(  \frac{\partial}{\partial t}-\Delta\right)  Z_{ab}  &  =2R_{cd}%
[D_{c}P_{dab}+D_{c}P_{dba}]\newline +2R_{acbd}M_{cd}\\
&  +2P_{acd}P_{bcd}-4P_{acd}P_{bdc}\newline +2R_{cd}R_{ce}R_{adbe}-{\frac
{1}{2t^{2}}}\,R_{ab}\\
&  +\left(
\begin{array}
[c]{c}%
\left(  S_{ijrs}S_{mntu}S_{klpq}+S_{ijmn}S_{klrs}S_{pqtu}\right)  \nabla
_{v}R_{rstu}\nabla_{v}R_{mnpq}\\
+2S_{ijmn}S_{klpq}(B_{mnpq}-B_{mnqp}+B_{mpnq}-B_{mqnp})
\end{array}
\right)  P_{ija}P_{klb}\\
&  -2S_{ijkl}P_{klb}\left(  -R_{de}D_{d}R_{ijae}+R_{idje}P_{dea}%
+R_{idae}P_{dje}\newline +R_{jdae}P_{ide}\right) \\
&  -2S_{ijkl}P_{ija}\left(  -R_{de}D_{d}R_{klbe}+R_{kdle}P_{deb}%
+R_{kdbe}P_{dle}\newline +R_{ldbe}P_{kde}\right) \\
&  +2S_{ijkl}\nabla_{p}P_{ija}\nabla_{p}P_{klb}-2S_{ijmn}S_{klpq}\nabla
_{r}R_{mnpq}\left(  \nabla_{r}P_{ija}P_{klb}+P_{ija}\nabla_{r}P_{klb}\right)
\end{align*}
where $B_{abcd}=R_{aebf}R_{cedf}\,.$ Taking into account symmetries and
combining some like terms, we have%
\begin{align*}
\left(  \frac{\partial}{\partial t}-\Delta\right)  Z_{ab}  &  =2R_{acbd}%
Z_{cd}\\
&  +2S_{ijkl}\nabla_{p}P_{ija}\nabla_{p}P_{klb}\\
&  -2S_{ijmn}S_{klpq}\nabla_{r}R_{mnpq}\left(  \nabla_{r}P_{ija}%
P_{klb}+P_{ija}\nabla_{r}P_{klb}\right) \\
&  +\left(  S_{ijrs}S_{mntu}S_{klpq}+S_{ijmn}S_{klrs}S_{pqtu}\right)
P_{ija}P_{klb}\nabla_{v}R_{rstu}\nabla_{v}R_{mnpq}\\
&  +2R_{acbd}S_{ijkl}P_{ijc}P_{kld}+2P_{acd}P_{bcd}-4P_{acd}P_{bdc}\\
&  +2S_{ijmn}S_{klpq}(B_{mnpq}-B_{mnqp}+B_{mpnq}-B_{mqnp})P_{ija}P_{klb}\\
&  -2\left(  S_{dekl}R_{diej}+S_{ijde}R_{dkel}\right)  P_{ija}\newline
P_{klb}\newline \\
&  -4S_{ijkl}R_{idae}P_{klb}P_{dje}\newline -4S_{ijkl}R_{kdbe}P_{ija}%
P_{dle}\newline \\
&  +2R_{cd}[D_{c}P_{dab}+D_{c}P_{dba}]\newline +2R_{cd}R_{ce}R_{adbe}%
-{\frac{1}{2t^{2}}}\,R_{ab}\\
&  +2S_{ijkl}P_{klb}R_{de}D_{d}R_{ijae}+2S_{ijkl}P_{ija}R_{de}D_{d}R_{klbe}.
\end{align*}
Separating the terms quadratic in $P_{abc}$ from the rest and noticing that
most of the rest is a square, we obtain:%
\begin{align*}
&  \left(  \frac{\partial}{\partial t}-\Delta\right)  Z_{ab}\\
&  =2R_{acbd}Z_{cd}\\
&  +2S_{mntu}\left(  S_{ijrs}P_{ija}\nabla_{v}R_{rstu}-\nabla_{v}%
P_{tua}+R_{tuaw}R_{vw}\right) \\
&  \times\left(  S_{klpq}P_{klb}\nabla_{v}R_{pqmn}-\nabla_{v}P_{mnb}%
+R_{mnbx}R_{vx}\right) \\
&  +2R_{acbd}S_{ijkl}P_{ijc}P_{kld}+2P_{acd}P_{bcd}-4P_{acd}P_{bdc}%
-4S_{ijkl}R_{idae}P_{klb}P_{dje}\newline \\
&  -4S_{ijkl}R_{kdbe}P_{ija}P_{dle}\newline +2S_{ijmn}S_{klpq}(B_{mnpq}%
-B_{mnqp}+B_{mpnq}-B_{mqnp})P_{ija}P_{klb}\\
&  -2\left(  S_{dekl}R_{diej}+S_{ijde}R_{dkel}\right)  P_{ija}\newline
P_{klb}\newline -{\frac{1}{2t^{2}}}\,R_{ab}.
\end{align*}
We simplify the terms that are quadratic in $P_{abc}\,$:%
\begin{align*}
A\doteqdot &  2R_{acbd}S_{ijkl}P_{ijc}P_{kld}+2P_{acd}P_{bcd}-4P_{acd}%
P_{bdc}-4S_{ijkl}R_{idae}P_{klb}P_{dje}\newline \\
&  +2S_{ijmn}S_{klpq}(B_{mnpq}-B_{mnqp}+B_{mpnq}-B_{mqnp})P_{ija}P_{klb}\\
&  -4S_{ijkl}R_{kdbe}P_{ija}P_{dle}-2\left(  S_{dekl}R_{diej}+S_{ijde}%
R_{dkel}\right)  P_{ija}\newline P_{klb}\newline \\
&  =2R_{acbd}S_{ijkl}P_{ijc}P_{kld}+2P_{acd}P_{bcd}-4P_{acd}P_{bdc}\\
&  -4S_{ijkl}R_{idae}P_{klb}P_{dje}\newline -4S_{ijkl}R_{kdbe}P_{ija}%
P_{dle}\newline +P_{cda}P_{cdb}\\
&  +2S_{ijmn}S_{klpq}(B_{mpnq}-B_{mqnp})P_{ija}P_{klb}-P_{kla}\newline
P_{klb}\newline -P_{ija}\newline P_{ijb}\\
&  =-2P_{acd}P_{bdc}+2R_{acbd}S_{ijkl}P_{ijc}P_{kld}-4S_{ijkl}R_{idae}%
P_{klb}P_{dje}\newline \\
&  -4S_{ijkl}R_{idbe}P_{kla}P_{dje}+2R_{iepf}R_{jeqf}S_{ijkl}P_{kla}%
S_{pqmn}P_{mnb}\\
&  -2R_{ieqf}R_{jepf}S_{ijkl}P_{kla}S_{pqmn}P_{mnb}\\
&  =-2P_{acd}P_{bdc}+2R_{acbd}S_{ijkl}P_{ijc}P_{kld}-4S_{ijkl}R_{idae}%
P_{klb}P_{dje}\newline \\
&  -4S_{ijkl}R_{idbe}P_{kla}P_{dje}+4R_{iepf}R_{jeqf}S_{ijkl}P_{kla}%
S_{pqmn}P_{mnb}.
\end{align*}
In the above computation we used the following identities. First, the
quadratic $B$ in curvature can be rewritten using%
\begin{align*}
R_{mnpq}^{2}  &  =R_{mnef}R_{pqef}=\left(  R_{menf}-R_{nemf}\right)  \left(
R_{peqf}-R_{qepf}\right) \\
&  =B_{mnpq}-B_{mnqp}-B_{nmpq}+B_{nmqp}\\
&  =2\left(  B_{mnpq}-B_{mnqp}\right)  .
\end{align*}
and
\begin{align*}
&  2S_{ijmn}S_{klpq}(B_{mpnq}-B_{mqnp})P_{ija}P_{klb}\\
&  =2S_{ijmn}S_{klpq}(R_{mepf}R_{neqf}-R_{meqf}R_{nepf})P_{ija}P_{klb}\\
&  =2S_{ijkl}S_{mnpq}R_{iepf}R_{jeqf}P_{kla}P_{mnb}-2S_{ijkl}S_{mnpq}%
R_{ieqf}R_{jepf}P_{kla}P_{mnb}.
\end{align*}
Second, we have the identities%
\begin{align*}
2S_{dekl}R_{diej}P_{ija}  &  =P_{kla}\\
2S_{ijde}R_{dkel}P_{klb}  &  =P_{ijb},
\end{align*}
where the first follows from%
\[
S_{dekl}\left(  R_{diej}-R_{djei}\right)  P_{ija}=S_{dekl}R_{deij}%
P_{ija}=P_{kla}%
\]
which is obtained from%
\[
S_{dekl}R_{diej}P_{ija}=S_{dekl}\left(  R_{deij}+R_{djei}\right)  P_{ija}%
\]
and the second is proved similarly. Third, we have the identity%
\begin{align*}
P_{cda}P_{cdb}  &  =\left(  P_{adc}+P_{cad}\right)  \left(  P_{bdc}%
+P_{cbd}\right) \\
&  =2\left(  P_{acd}P_{bcd}-P_{acd}P_{bdc}\right)  .
\end{align*}

Now define $E_{ija}\doteqdot S_{ijkl}P_{kla}$ to simplify notation further.
Then%
\begin{align*}
A  &  =-2P_{acd}P_{bdc}+2R_{acbd}S_{ijkl}P_{ijc}P_{kld}-4R_{idae}%
E_{ijb}P_{dje}\newline \\
&  -4R_{idbe}E_{ija}P_{dje}+4R_{iepf}R_{jeqf}E_{ija}E_{pqb}\,.
\end{align*}
Recall that%
\begin{align*}
R_{abcd}  &  =\sum_{N}Y_{ab}^{N}Y_{cd}^{N}\\
P_{abc}  &  =\sum_{N}Y_{ab}^{N}X_{c}^{N}\\
M_{ab}  &  =\sum_{N}X_{a}^{N}X_{b}^{N}.
\end{align*}
Thus%
\begin{align*}
A  &  =-2P_{acd}P_{bdc}+2R_{acbd}S_{ijkl}P_{ijc}P_{kld}+4Y_{pf}^{N}Y_{ah}%
^{N}Y_{qf}^{M}X_{h}^{M}E_{pqb}\newline \\
&  +4Y_{id}^{N}Y_{be}^{N}Y_{jd}^{M}X_{e}^{M}E_{ija}+4Y_{id}^{N}Y_{pf}%
^{N}Y_{jd}^{M}Y_{qf}^{M}E_{ija}E_{pqb}\\
&  =4\left(  Y_{id}^{N}Y_{jd}^{M}E_{ija}+Y_{ah}^{N}X_{h}^{M}\newline \right)
\left(  Y_{pf}^{N}Y_{qf}^{M}E_{pqb}+Y_{be}^{N}X_{e}^{M}\right)  -4Y_{ah}%
^{N}X_{h}^{M}Y_{be}^{N}X_{e}^{M}\\
&  -2P_{acd}P_{bdc}+2R_{acbd}S_{ijkl}P_{ijc}P_{kld}\\
&  =4\left(  Y_{id}^{N}Y_{jd}^{M}E_{ija}+Y_{ah}^{N}X_{h}^{M}\newline \right)
\left(  Y_{pf}^{N}Y_{qf}^{M}E_{pqb}+Y_{be}^{N}X_{e}^{M}\right)  -2R_{acbd}%
Z_{cd}\\
&  -2R_{acbd}M_{cd}-2P_{acd}P_{bdc}.
\end{align*}
Note that the terms $Y_{id}^{N}Y_{jd}^{M}E_{ija}$ and $Y_{pf}^{N}Y_{qf}%
^{M}E_{pqb}$ are both antisymmetric in $N$ and $M.$

Now the terms on the last line may be rewritten as%
\begin{align*}
&  -2R_{acbd}M_{cd}-2P_{acd}P_{bdc}\\
&  =-2Y_{ac}^{N}Y_{bd}^{N}X_{c}^{M}X_{d}^{M}-2Y_{ac}^{N}X_{d}^{N}Y_{bd}%
^{M}X_{c}^{M}P_{acd}P_{bdc}\\
&  =-\left(  Y_{ac}^{N}X_{c}^{M}+Y_{ac}^{M}X_{c}^{N}\right)  \left(
Y_{bd}^{N}X_{d}^{M}+Y_{bd}^{M}X_{d}^{N}\right)
\end{align*}
so that
\begin{align*}
A  &  =4\left(  Y_{id}^{N}Y_{jd}^{M}E_{ija}+Y_{ah}^{N}X_{h}^{M}\newline
\right)  \left(  Y_{pf}^{N}Y_{qf}^{M}E_{pqb}+Y_{be}^{N}X_{e}^{M}\right)
-2R_{acbd}Z_{cd}\\
&  -\left(  Y_{ac}^{N}X_{c}^{M}+Y_{ac}^{M}X_{c}^{N}\right)  \left(  Y_{bd}%
^{N}X_{d}^{M}+Y_{bd}^{M}X_{d}^{N}\right) \\
&  =\left(  2Y_{id}^{N}Y_{jd}^{M}E_{ija}+Y_{ah}^{N}X_{h}^{M}\newline
-Y_{ah}^{M}X_{h}^{N}\newline \right)  \left(  2Y_{pf}^{N}Y_{qf}^{M}%
E_{pqb}+Y_{be}^{N}X_{e}^{M}-Y_{be}^{M}X_{e}^{N}\right) \\
&  -2R_{acbd}Z_{cd},
\end{align*}
where we used
\[
Y_{ah}^{N}X_{h}^{M}\left(  Y_{be}^{N}X_{e}^{M}+Y_{be}^{M}X_{e}^{N}\right)
=\frac{1}{2}\left(  Y_{ah}^{N}X_{h}^{M}\newline +Y_{ah}^{M}X_{h}^{N}%
\newline \right)  \left(  Y_{be}^{N}X_{e}^{M}+Y_{be}^{M}X_{e}^{N}\right)
\]
and%
\[
Y_{id}^{N}Y_{jd}^{M}E_{ija}\left(  Y_{be}^{N}X_{e}^{M}+Y_{be}^{M}X_{e}%
^{N}\right)  =0
\]
(since $Y_{be}^{N}X_{e}^{M}+Y_{be}^{M}X_{e}^{N}$ is symmetric in $N$ and $M$).

Therefore the evolution equation for $Z_{ab}$ is:%
\begin{align*}
&  \left(  \frac{\partial}{\partial t}-\Delta\right)  Z_{ab}\\
&  =2S_{mntu}\left(  S_{ijrs}P_{ija}\nabla_{v}R_{rstu}-\nabla_{v}%
P_{tua}+R_{tuaw}R_{vw}\right) \\
&  \times\left(  S_{klpq}P_{klb}\nabla_{v}R_{pqmn}-\nabla_{v}P_{mnb}%
+R_{mnbx}R_{vx}\right) \\
&  +\left(  2Y_{id}^{N}Y_{jd}^{M}E_{ija}+Y_{ah}^{N}X_{h}^{M}\newline
-Y_{ah}^{M}X_{h}^{N}\newline \right)  \left(  2Y_{pf}^{N}Y_{qf}^{M}%
E_{pqb}+Y_{be}^{N}X_{e}^{M}-Y_{be}^{M}X_{e}^{N}\right) \\
&  -{\frac{1}{2t^{2}}}\,R_{ab},
\end{align*}
where $E_{ija}=S_{ijkl}P_{kla}.$ We may rewrite this as%
\begin{align*}
&  \left(  \frac{\partial}{\partial t}-\Delta\right)  Z_{ab}\\
&  =2S_{mntu}\left(  S_{ijrs}P_{ija}\nabla_{v}R_{rstu}-\nabla_{v}%
P_{tua}+R_{tuaw}R_{vw}+\frac{1}{2t}R_{tuav}\right) \\
&  \times\left(  S_{klpq}P_{klb}\nabla_{v}R_{pqmn}-\nabla_{v}P_{mnb}%
+R_{mnbx}R_{vx}+\frac{1}{2t}R_{mnbv}\right) \\
&  +\left(  2Y_{id}^{N}Y_{jd}^{M}E_{ija}+Y_{ah}^{N}X_{h}^{M}\newline
-Y_{ah}^{M}X_{h}^{N}\newline \right)  \left(  2Y_{pf}^{N}Y_{qf}^{M}%
E_{pqb}+Y_{be}^{N}X_{e}^{M}-Y_{be}^{M}X_{e}^{N}\right) \\
&  -\frac{1}{t}\left(  S_{ijrs}P_{ija}\nabla_{v}R_{rsbv}-\nabla_{v}%
P_{bva}+R_{bvaw}R_{vw}+\frac{1}{2t}R_{bvav}\right) \\
&  -\frac{1}{t}\left(  S_{klpq}P_{klb}\nabla_{v}R_{pqav}-\nabla_{v}%
P_{avb}+R_{avbx}R_{vx}+\frac{1}{2t}R_{avbv}\right)
\end{align*}

Now%
\begin{align*}
-\nabla_{v}P_{bva}+R_{bvaw}R_{vw}+\frac{1}{2t}R_{bvav} &  =-\nabla_{v}\left(
\nabla_{b}R_{va}-\nabla_{v}R_{ba}\right)  +R_{bvaw}R_{vw}+\frac{1}{2t}R_{ba}\\
&  =\Delta R_{ba}-\nabla_{b}\nabla_{v}R_{va}-R_{bw}R_{wa}+R_{bvaw}R_{vw}\\
&  +R_{bvaw}R_{vw}+\frac{1}{2t}R_{ba}\\
&  =\Delta R_{ba}-\frac{1}{2}\nabla_{b}\nabla_{a}R+2R_{bvaw}R_{vw}%
-R_{ac}R_{bc}+\frac{1}{2t}R_{ab}\\
&  =M_{ab}%
\end{align*}
and%
\[
\nabla_{v}R_{rsbv}=\nabla_{r}R_{vsbv}+\nabla_{s}R_{rvbv}=-\nabla_{r}%
R_{sb}+\nabla_{s}R_{rb}=-P_{rsb}.
\]
Hence%
\begin{align*}
&  \left(  \frac{\partial}{\partial t}-\Delta\right)  Z_{ab}\\
&  =2S_{mntu}\left(  S_{ijrs}P_{ija}\nabla_{v}R_{rstu}-\nabla_{v}%
P_{tua}+R_{tuaw}R_{vw}+\frac{1}{2t}R_{tuav}\right)  \\
&  \times\left(  S_{klpq}P_{klb}\nabla_{v}R_{pqmn}-\nabla_{v}P_{mnb}%
+R_{mnbx}R_{vx}+\frac{1}{2t}R_{mnbv}\right)  \\
&  +\left(  2Y_{id}^{N}Y_{jd}^{M}E_{ija}+Y_{ah}^{N}X_{h}^{M}\newline
-Y_{ah}^{M}X_{h}^{N}\newline \right)  \left(  2Y_{pf}^{N}Y_{qf}^{M}%
E_{pqb}+Y_{be}^{N}X_{e}^{M}-Y_{be}^{M}X_{e}^{N}\right)  \\
&  -\frac{1}{t}\left(  -S_{ijrs}P_{ija}P_{rsb}+M_{ab}\right)  -\frac{1}%
{t}\left(  -S_{klpq}P_{klb}P_{pqa}+M_{ab}\right)  .
\end{align*}
That is
\begin{align*}
&  \left(  \frac{\partial}{\partial t}-\Delta\right)  Z_{ab}=2S_{mntu}\left(
S_{ijrs}P_{ija}\nabla_{v}R_{rstu}-\nabla_{v}P_{tua}+R_{tuaw}R_{vw}+\frac
{1}{2t}R_{tuav}\right)  \\
&  \times\left(  S_{klpq}P_{klb}\nabla_{v}R_{pqmn}-\nabla_{v}P_{mnb}%
+R_{mnbx}R_{vx}+\frac{1}{2t}R_{mnbv}\right)  \\
&  +\left(  2Y_{id}^{N}Y_{jd}^{M}E_{ija}+Y_{ah}^{N}X_{h}^{M}\newline
-Y_{ah}^{M}X_{h}^{N}\newline \right)  \left(  2Y_{pf}^{N}Y_{qf}^{M}%
E_{pqb}+Y_{be}^{N}X_{e}^{M}-Y_{be}^{M}X_{e}^{N}\right)  \\
&  -\frac{2}{t}Z_{ab}.
\end{align*}
This completes the proof of the main theorem.\medskip

Now it is interesting to observe the form of the resulting evolution equation
when we trace the above computation. Define%
\[
Z\doteqdot Z_{aa}=\frac{1}{2}\left(  \Delta R+2\left|  Rc\right|  ^{2}%
+\frac{1}{t}R\right)  -S_{ijkl}P_{ija}P_{kla}.
\]
This trace Harnack quantity is slightly different than
(\ref{trace-inequality1}); however its positivity is of course still a
corollary of Hamilton's matrix Harnack estimate. Tracing equation (\ref{MT}),
we have%
\[
\left(  \frac{\partial}{\partial t}-\Delta\right)  Z=2S_{ijkl}K_{mnij}%
K_{mnkl}+\left(  L_{a}^{NM}\right)  ^{2}-\frac{2}{t}Z\geq-\frac{2}{t}Z.
\]
Note however that an application of the maximum principle does not yield a
\emph{direct} proof of the trace Harnack inequality. This is because in order
for $L_{a}^{NM},$ which appears on the rhs, to be well-defined, we need the
matrix Harnack quantity to be nonnegative!\medskip

\textbf{Acknowledgment. }We would like to express our gratitude to Richard
Hamilton for suggesting to us the method of proof given in this paper.

\end{document}